\documentclass[11pt,leqno,twoside]{article}

\usepackage{amsfonts,amsmath,amsthm,amssymb}
\usepackage{enumerate}
\usepackage{graphics,graphicx,subfigure}

\usepackage{color}
\usepackage{todonotes}
\usepackage{cancel}
\usepackage{url}
\usepackage{hyperref}
\usepackage{makeidx}
\usepackage{showidx}
\usepackage{multicol}        % used for the two-column index
\usepackage{xspace}
\usepackage{stmaryrd}        % for \llbracket, \rrbracket
\usepackage{pifont}          % for \ding{227}
\usepackage{fancybox}        % for oval boxes
\usepackage{bm}

 \usepackage{fancyhdr}
\theoremstyle{plain}
 \theoremstyle{definition}
 \newtheorem{lem}{Lemma}
 \newtheorem{defn}[lem]{Definition}
 \newtheorem{thm}[lem]{Theorem}
 \newtheorem{prop}[lem]{Proposition}
 \newtheorem{cor}[lem]{Corollary}
 \newtheorem{notn}[lem]{Notations}
 \newtheorem{pb}[lem]{Problem}
 \newtheorem{form}[lem]{Formulae}
 
 \newtheorem*{rk}{Remark}
 \newtheorem*{com}{Comment}
 \newtheorem*{ex}{Example}
 \theoremstyle{remark}

 \newcommand{\blem}{\begin{lem}}
 \newcommand{\elem}{\end{lem}}
 \newcommand{\bdefn}{\begin{defn}}
 \newcommand{\edefn}{\end{defn}}
 \newcommand{\bthm}{\begin{thm} }
 \newcommand{\ethm}{\end{thm}}
 \newcommand{\bprop}{\begin{prop}}
 \newcommand{\eprop}{\end{prop}}
 \newcommand{\bcor}{\begin{cor}}
 \newcommand{\ecor}{\end{cor}}
 \newcommand{\bnotn}{\begin{notn}}
 \newcommand{\enotn}{\end{notn}}
 \newcommand{\bpb}{\begin{pb}}
 \newcommand{\epb}{\end{pb}}
 \newcommand{\bform}{\begin{form}}
 \newcommand{\eform}{\end{form}}
 \newcommand{\brk}{\begin{rk}}
 \newcommand{\erk}{\end{rk}}
 \newcommand{\bcom}{\begin{com}}
 \newcommand{\ecom}{\end{com}}
 \newcommand{\bex}{\begin{ex}}
 \newcommand{\eex}{\end{ex}}
 \newcommand{\bpf}{\begin{proof}}
 \newcommand{\epf}{\end{proof}}

\newcommand{\cC}{\mathcal{C}}

\newcommand{\cP}{\mathcal{P}}

\newcommand{\cU}{\mathcal{U}}
\newcommand{\cV}{\mathcal{V}}

\newcommand{\bN}{\mathbb{N}}

\newcommand{\bR}{\mathbb{R}}

\newcommand{\bT}{\mathbb{T}}

\newcommand{\bZ}{\mathbb{Z}}

\newcommand{\be}{\begin{equation}}
\newcommand{\ee}{\end{equation}}
\newcommand{\bal}{\begin{align}}
\newcommand{\eal}{\end{align}}
\newcommand{\ba}{\begin{align*}}
\newcommand{\ea}{\end{align*}}
\newcommand{\bmx}{\begin{matrix}}
\newcommand{\emx}{\end{matrix}}
\newcommand{\bbmx}{\begin{bmatrix}}
\newcommand{\ebmx}{\end{bmatrix}}
\newcommand{\bpmx}{\begin{pmatrix}}
\newcommand{\epmx}{\end{pmatrix}}
\newcommand{\bvmx}{\begin{vmatrix}}
\newcommand{\evmx}{\end{vmatrix}}

\newcommand{\wh}{\widehat}
\newcommand{\wt}{\widetilde}
\newcommand{\f}{\frac}
\newcommand{\df}{\dfrac}

\newcommand{\imp}{\Longrightarrow}

\newcommand{\inc}{\subseteq}

\newcommand{\lbt}{\llbracket}
\newcommand{\rbt}{\rrbracket}

\newcommand{\la}{\lambda}

\newcommand{\eps}{\varepsilon}
  
\newcommand{\ibt}[2]{\lbt {#1} \hspace{-1mm}: \hspace{-1mm} {#2} \rbt }

\title{\vspace{-20mm}Determining  Projection Constants of Univariate Polynomial Spaces\medskip\hrule height 1.2pt \vspace{-6mm}}
\author{Simon Foucart\footnote{The research of the first author is partially funded by the NSF grant DMS-1622134.} \, (Texas A\&M University) and Jean B. Lasserre\footnote{The research of the second author is funded by the European Research Council (ERC) under
the European's Union Horizon 2020 research and innovation program (grant agreement 666981 TAMING).} \, (University of Toulouse)}
\date{\vspace{-6mm}\rule{100mm}{0.8pt}}

\newcommand\shorttitle{Determining the projection constants of univariate polynomial spaces}
\newcommand\authors{S. Foucart, J. B. Lasserre}

\begin{document}
\maketitle

%% Add abstract, keywords, and AMS classification
%\vspace{-15mm}
\begin{abstract}
The long-standing problem of minimal projections is addressed from a computational point of view.
Techniques to determine bounds on the projection constants of univariate polynomial spaces are presented.
The upper bound, produced by a linear program,
and the lower bound, produced by a semidefinite program exploiting the method of moments,
are often close enough to deduce the projection constant with reasonable accuracy.
The implementation of these programs makes it possible to find the projection constant of several three-dimensional spaces with five digits of accuracy,
as well as the projection constants of the spaces of cubic, quartic, and quintic polynomials with four digits of accuracy.
Beliefs about uniqueness and shape-preservation of minimal projections are contested along the way.
\end{abstract}



\vspace{-5mm}
\begin{center}
\rule{100mm}{0.8pt}
\end{center}

%%%%%%%%%%%%%%%%%
%% The main text starts here %%
%%%%%%%%%%%%%%%%%

\section{Introduction}

The problem of minimal projections has attracted the attention of approximation theorists for about half a century.
The survey of Cheney and Price \cite{ChePri} still provides a fine account on the topic.
The fact that the Fourier projection is uniquely minimal from $\cC(\bT)$ onto the space of trigonometric polynomials of degree at most $d$,
derived from Berman--Marcinkiewicz formula,
undoubtedly stands as a highlight of the subject, see \cite{Loz,CHMSW}.
But when the focus is put on algebraic rather than trigonometric polynomials,
the situation becomes dramatically more complicated.
Besides the trivial cases of degree $d=0$ and $d=1$,
only the case $d=2$ has been resolved,
albeit at the cost of considerable efforts
deployed by Chalmers and Metcalf \cite{ChaMet}.
In fact, traditional analyses may have reached their limitation for the problem of minimal projections.
The present article is our way of advocating a change of philosophy 
to promote the integration in classical Approximation Theory 
of modern optimization methods,
in particular methods based on moments and positive polynomials,
see \cite{L}.
Such techniques have already been used for minimal projections \cite{NFAO} and also in the context of constrained approximation \cite{Basc}.

Before presenting the results and insights generated by our approach,
let us formalize the problem of interest precisely.
Throughout the article, we consider an $M$-dimensional subspace $\cU$ of the space $\cC[-1,1]$ of continuous functions defined on the interval $[-1,1]$.
Typically, the space $\cU$ consists of algebraic polynomials.
A projection $P$ from $\cC[-1,1]$ onto $\cU$ is just a linear map from $\cC[-1,1]$ onto~$\cU$
that reproduces $\cU$,
i.e., that satisfies $P(u) = u$ for all $u \in \cU$.
A minimal projection from $\cC[-1,1]$ onto $\cU$ is a projection $P$ from $\cC[-1,1]$ onto $\cU$ with minimal norm
\be
\label{1}
\|P\|_{\infty \to \infty}
= \max_{ \|f\|_\infty \le 1}  \|P(f)\|_\infty 
= \max_{x \in [-1,1]} \max_{\|f\|_\infty \le 1} | P(f)(x)|.
\ee
The projection constant $\la(\cU)$ of $\cU$ (relative to $\cC[-1,1]$) is the value of the minimum, i.e.,
\be
\label{Ori1}
\la(\cU) :=
\min_{P} \; \|P\|_{\infty \to \infty}
\quad
\mbox{s.to $P$ being a projection from $\cC[-1,1]$ onto $\cU$}.
\ee
Fixing a basis $(u_1,\ldots,u_M)$ for $\cU$,
any projection $P$ from $\cC[-1,1]$ onto $\cU$ can be represented 
by some linear functionals $\eta_1,\ldots,\eta_M$ on $\cC[-1,1]$
as
\be
P(f) = \sum_{m=1}^M \eta_m(f) u_m,
\qquad 
\mbox{with}
\quad
\eta_m(u_{m'}) = \delta_{m,m'}, \quad
m,m' \in \ibt{1}{M} .
\ee
Moreover, any bounded linear functional $\eta$ on $\cC[-1,1]$ can be represented 
by some signed Borel measure $\mu$ on $[-1,1]$ as
\be
\eta(f) = \int_{-1}^1 f d\mu.
\ee
It follows easily that the projection constant of $\cU$ can also be written as 
\be
\label{Original}
\la(\cU) =  \inf_{\mu_1,\ldots,\mu_M}
\max_{x \in [-1,1]} \int_{-1}^1 \left| \sum_{m=1}^M u_m(x) d\mu_m \right|
\quad 
\mbox{s.to }
\int_{-1}^1 u_{m'} d\mu_{m} = \delta_{m,m'},
\; m,m' \in \ibt{1}{M},
\ee
where the infimum (in fact,  minimum)  is taken over all signed Borel measures $\mu_1,\ldots,\mu_M$ on $[-1,1]$.
There is of course a difficulty originating from the infinite-dimensionality of the optimization variables $\mu_1,\ldots,\mu_M$,
so the optimization program \eqref{Original} cannot (a priori) be performed exactly.

Our strategy  consists in producing computable upper bounds 
(see Section \ref{SecUB})
and lower bounds (see Section~\ref{SecLB})
that are sufficiently close to determine the value of the projection constant with, say, four or five digits of accuracy
(corresponding here to three or four digits after the decimal point).
The computational burden can be lighten by exploiting some symmetry properties of minimal projections (see Section~\ref{SecSym}).
Implementing our strategy enables us, for instance, to:\vspace{-5mm}
\begin{itemize}
\item retrieve numerically the result of \cite{ChaMet} for quadratic polynomials, namely
\be
\la(\cP_2) \approx 1.2201,
\ee
and allude\footnote{\label{FN}We are not being more assertive here because, strictly speaking, our computed minimal projection is not an exact minimal projection.} to the nonuniqueness of minimal projections onto the quadratics;
\item determine with five digits of accuracy the projection constants of several other three-dimensional polynomial spaces, and note in passing that
\be 
\la({\rm span}\{1,x^2,x^3\}) = 1;
\ee
\item determine with four digits of accuracy the projection constants of the spaces of cubic, quatric, and quintic polynomials, namely
\begin{align} 
\la(\cP_3) &\approx 1.365,\\
\la(\cP_4) & \approx 1.459,\\
\la(\cP_5) & \approx 1.538,
\end{align}
and hint\footnote{Same reservation as in footnote \ref{FN}.} that minimal projections onto $\cP_d$ do not in general preserve $d$-convexity,
thus \mbox{disproving} a conjecture from \cite{PCM}; 
\item provide state-of-the-art upper and lower bounds for the projection constants of the spaces of polynomials of degree at most $d$ until $d = 12$.
\end{itemize}\vspace{-5mm}
All of our results can be reproduced by downloading the {\sc matlab} code available on the authors' webpages.
The packages {\sf CVX} and {\sf Chebfun} are required to execute the code.

\section{Computable upper bound}
\label{SecUB}

We present in this section a discretization of the problem \eqref{Original} that leads to a computable upper bound for the projection constant $\la(\cU)$ of a polynomial subspace $\cU$ of $\cC[-1,1]$.
This involves in fact two discretizations:
one for the signed Borel measures $\mu_1,\ldots,\mu_M$
and one for the interval $[-1,1]$.
We start by discretizing the measures.

\bprop
\label{PropUB1}
Given a basis $(u_1,\ldots,u_M)$ for a subspace $\cU$ of $\cC[-1,1]$,
if $v_1,\ldots,v_K$ are distinct points in $[-1,1]$,
then
\be
\label{UB1}
\la(\cU) \le \min_{A \in \bR^{M \times K}} 
\max_{x \in [-1,1]} 
\sum_{k=1}^K \left| \sum_{m=1}^M A_{m,k} u_m(x) \right|
\;
\mbox{s.to }
 \sum_{k=1}^K A_{m,k} u_{m'}(v_k) = \delta_{m,m'},
 \,  m,m' \in \ibt{1}{M} .
\ee
\eprop

\bpf
We obtain an upper bound for $\la(\cU)$
whenever the minimization in \eqref{Original} is carried over a subset of all signed Borel measures $\mu_1,\ldots,\mu_M$.
In particular, we choose here the subset of all linear combinations of Dirac measures at $v_1,\ldots,v_K$,
i.e., measures  of the form
\be
\mu_m = \sum_{k=1}^K A_{m,k} \delta_{v_k}, 
\qquad m \in \ibt{1}{M}.
\ee
Under this restriction, the duality constraints in \eqref{Original} read,
for all $m,m' \in \ibt{1}{M}$,
\be
\label{Exp1}
\delta_{m,m'} = 
\int_{-1}^1 u_{m'} d\left( \sum_{k=1}^K A_{m,k} \delta_{v_k} \right)
=\sum_{k=1}^K A_{m,k} \int_{-1}^1 u_{m'}  d\delta_{v_k}
= \sum_{k=1}^K A_{m,k} u_{m'}(v_k),
\ee
while the integral appearing in the objective function becomes,
 for a fixed $x \in [-1,1]$, 
\begin{align}
\label{Exp2}
\int_{-1}^1 \left| \sum_{m=1}^M u_m(x) d \left( \sum_{k=1}^K A_{m,k} \delta_{v_k} \right) \right|
& = \int_{-1}^1 \left| \sum_{k=1}^K  \left( \sum_{m=1}^M A_{m,k} u_m(x) \right) d\delta_{v_k} \right| \\
\nonumber
& = \int_{-1}^{1} \sum_{k=1}^K \left| \sum_{m=1}^M A_{m,k} u_m(x) \right| d\delta_{v_k}
= \sum_{k=1}^K \left| \sum_{m=1}^M A_{m,k} u_m(x) \right|.
\end{align}
Taking the expressions \eqref{Exp1} and \eqref{Exp2} into account yields the upper bound \eqref{UB1}.
\epf

Proposition \ref{PropUB1} is not directly exploitable
due to the presence of the maximum over the interval $[-1,1]$.
We can replace it by a maximum over a discretized grid,
as long as we are able to bound the maximum over $[-1,1]$ by the maximum over this grid.
For polynomials,
the comparison between the discrete and continuous max-norms is a well-studied topic,
especially for equispaced points (see e.g.
\cite{EZ,CR,Rak}).
But equispaced points are not the most suitable,
since the two norms are comparable when the number of points scales quadratically with the degree.
In contrast,
for zeros of Chebyshev polynomials,
the number of points only needs to scale linearly with the degree.
Here is a quantitative version of this assertion,
which can be found in \cite{EZ}.

\blem
\label{LemDiscVsCont}
Let $w_1> \cdots > w_L \in [-1,1]$ be the Chebyshev zeros given by $w_\ell = \cos(\theta_\ell)$, where $\theta_\ell = \pi(\ell-1/2)/L$.
For any algebraic polynomial $p$ of degree at most $d$, one has
\be
\label{DiscVsCont}
\max_{x \in [-1,1]}|p(x)| \le \cos\left(\df{\pi}{2} \df{d}{L} \right)^{-1} \max_{\ell \in \lbt 1:L \rbt} |p(w_\ell)|.
\ee
\elem

We are now in a position to derive the awaited computable upper bound for the projection constant.

\bprop
\label{PropUB}
Given a basis $(u_1,\ldots,u_M)$ for a subspace $\cU$ of $\cC[-1,1]$ consisting of polynomials of degree at most $d$,
let $v_1,\ldots,v_K$ be distinct points in $[-1,1]$.
With $ w_1 > \cdots > w_L$ denoting the Chebyshev zeros
 $w_\ell = \cos(\pi (\ell-1/2)/L)$ and with $\rho = \cos\left( (\pi d)/(2 L) \right)^{-1} \ge 1$, one has
\be
\label{DD}
\la(\cU)
\le \rho \times
\hspace{-1.5mm} \min_{A \in \bR^{M \times K}} 
\max_{\ell \in \lbt 1 : L \rbt } 
\sum_{k=1}^K \left| \sum_{m=1}^M \hspace{-.5mm}A_{m,k} u_m(w_\ell) \right|
\;
\mbox{s.to }
 \sum_{k=1}^K \hspace{-.5mm}A_{m,k} u_{m'}(v_k) = \delta_{m,m'},
\,   m,m' \hspace{-.5mm} \in \hspace{-.5mm} \ibt{1}{M}.
\ee 
\eprop

\bpf
According to Proposition \ref{PropUB1}, it suffices to show that, for any $A \in \bR^{M \times K}$,
\be
\label{ObjUB2}
\max_{x \in [-1,1]} 
\sum_{k=1}^K \left| \sum_{m=1}^M A_{m,k} u_m(x) \right|
\le \rho 
\max_{\ell \in \lbt 1: L \rbt} 
\sum_{k=1}^K \left| \sum_{m=1}^M A_{m,k} u_m(w_\ell) \right|.
\ee
For a fixed $x \in [-1,1]$, 
we can find signs $\eps_1,\ldots,\eps_K \in \{\pm 1\}$ such that
\begin{align}
\label{CompAbs}
\sum_{k=1}^K \left| \sum_{m=1}^M A_{m,k} u_m(x) \right|
& = \left| \sum_{k=1}^K \eps_k \sum_{m=1}^M A_{m,k} u_m(x)  \right|
 \le \rho 
\max_{\ell \in \lbt 1: L \rbt} 
\left| \sum_{k=1}^K \eps_k \sum_{m=1}^M A_{m,k} u_m(w_\ell) \right|\\
\nonumber & \le \rho 
\max_{\ell \in \lbt 1:L \rbt} 
\sum_{k=1}^K \left|  \sum_{m=1}^M A_{m,k} u_m(w_\ell) \right|,
\end{align}
where Lemma \ref{LemDiscVsCont} was used for the first inequality in \eqref{CompAbs}.
Taking the maximum over $x \in [-1,1]$
yields the desired inequality \eqref{ObjUB2}.
\epf

We close this section by highlighting how the upper bound from Proposition \ref{PropUB}
is effectively computed by solving a linear program.
For this purpose, we introduce the collocation matrices $V \in \bR^{K \times M}$
and $W \in \bR^{L \times M}$
of the basis $(u_1,\ldots,u_M)$
at the points $v_1,\ldots,v_K$
(usually chosen as equispaced points in $[-1,1]$) 
and at the Chebyshev zeros $w_1,\ldots,w_L$.
These matrices are defined by
\be
\label{DefVW1}
V_{k,m} = u_m(v_k)
\qquad \mbox{and} \qquad
W_{\ell,m} = u_m(w_{\ell}),
\qquad
k \in \ibt{1}{K}, \;
\ell \in \ibt{1}{L}, \;
m \in \ibt{1}{M}.
\ee
With this notation, the objective function and the constraints in \eqref{DD} read, respectively,
\be
\max_{\ell \in \lbt 1:L \rbt} \sum_{k=1}^K |(WA)_{\ell,k}|
\qquad \mbox{and} \qquad
(AV)_{m,m'} = \delta_{m,m'},
\quad m,m' \in \ibt{1}{M}.
\ee
Then, introducing slack variables through $B \in \bR^{L \times K}$
(such that $|(WA)_{\ell,k}| \le B_{\ell,k}$ for all $\ell,k$)
and $c \in \bR$
(such that $\sum_{k=1}^K B_{\ell,k} \le c$ for all $\ell$),
we arrive at the following linear optimization problem,
written in an easily implementable form.

\ovalbox{
\begin{minipage}{0.97\textwidth}
\medskip
\begin{center}
\textbf{Upper bound for the projection constant}\vspace{-2mm}\\
%\rule{4.6in}{.8pt}
\rule{0.96\textwidth}{.8pt}
\end{center}
Inputs: basis for a space $\cU$ of polynomials of degree $\le d$, parameters $K$ and $L$. 
\be
\label{ProgUB}
\la(\cU) 
\le \cos \left( \df{\pi}{2}\df{d}{L} \right)^{-1}\times
\min_{\substack{A \in \bR^{M \times K}\\ B \in \bR^{L \times K}\\ c \in \bR }}
\; c
 \qquad \mbox{s.to \;\;} 
\left\{  
\bmx
  AV = I_M,  \hfill \\
-B \le WA \le B,  \hfill \\
 B {\bf 1} \le c {\bf 1}, \hfill
\emx
\right.
\ee
where the matrices $V$ and $W$ are defined in \eqref{DefVW1}.
\medskip
\end{minipage}
}

\section{Computable lower bound}
\label{SecLB}

We present in this section a computable lower bound for the projection constant of a \mbox{polynomial} \mbox{subspace} $\cU$ of $\cC[-1,1]$.
It again involves a discretization of the problem \eqref{Original},
followed by an \mbox{application} of the moment method.
This method, based on classical moment problems
(see e.g.~\cite{Schm}),
characterizes a measure by its sequence of moments,
typically via a semidefinite condition.
The \mbox{discrete} Hamburger moment problem, for instance,
states that,
for a sequence $(y_k)_{k \ge 0}$ of real \mbox{numbers},
there exists a nonnegative measure $\mu$ on $\bR$ such that
\be
\int_{-\infty}^\infty x^k d\mu(x) = y_k,
\qquad
k \ge 0,
\ee
if and only if an infinite Hankel matrix is positive semidefinite, precisely 
\be 
{\rm Hank}_\infty(y) := \bbmx
y_0 & y_1 & y_2 & y_3 & \cdots\\
y_1 & y_2 & y_3 & \reflectbox{$\ddots$} & \\
y_2 & y_3 &\reflectbox{$\ddots$} & & \\
y_3 & \reflectbox{$\ddots$} & & & \\
\vdots & & & &
\ebmx
\succeq 0.
\ee
We could directly use this characterization to substitute the measures $\mu_1,\ldots,\mu_M$ by their sequences $y_1,\ldots,y_M$ of moments as optimization variables,
but we would need to add extra semidefinite conditions ensuring that the measures $\mu_1,\ldots,\mu_M$ are localized on $[-1,1]$.
Instead, we prefer to rely on the discrete trigonometric moment problem, which states\footnote{The classical statement, found e.g. in \cite[Theorem 11.3]{Schm},
concerns sequences of complex numbers indexed by $n \in \bZ$
and obtained as $\int_{|z|=1} z^{-n} d\mu(z)$ for some Radon measure on the unit circle.
We omit the verification that it implies the statement being made here.} that,
for a sequence $(y_k)_{k \ge 0}$ of real numbers,
there exists a nonnegative measure $\mu$ on $[0,\pi]$ such that
\be 
\int_{0}^\pi \cos( k \theta) d \mu(\theta)
 = y_k,
 \qquad 
 k \ge 0,
 \ee
 if and only if an infinite Toeplitz matix is positive semidefinite, precisely
 \be 
 {\rm Toep}_\infty(y) := \bbmx
y_0 & y_1 & y_2 & y_3 & \cdots\\
y_1 & y_0 & y_1 & y_2 & \\
y_1 & y_1 & y_0 & y_1 & \ddots \\
y_3 & y_2 & y_1 & \ddots & \ddots \\
\vdots & & \ddots & \ddots & \ddots
\ebmx
 \succeq 0.
 \ee
 There are two reasons explaining our preference:
 firstly, localization conditions are not necessary;
 secondly, the Toeplitz structure is numerically more favorable than the Hankel structure.

 Now, in order to invoke the discrete trigonometric moment problem,
 we first have to transform the expression of the projection constant given in \eqref{Original}. 
 This is done below.

\blem
\label{LemTransformed}
Given a basis $(u_1,\ldots,u_M)$ for a subspace $\cU$ of $\cC[-1,1]$, 
if the functions $\wh{u}_1,  \ldots, \wh{u}_M$ are defined on $[0,\pi]$ by
$\wh{u}_m(\theta) = u_m( \cos(\theta))$,
then
\be
\label{Transformed}
\la(\cU) =  \inf_{\mu_1,\ldots,\mu_M}
\max_{\theta \in [0,\pi]} \int_{0}^\pi \left| \sum_{m=1}^M \wh{u}_m(\theta) d\mu_m \right|
\quad 
\mbox{s.to }
\int_{0}^\pi \wh{u}_{m'} d\mu_{m} = \delta_{m,m'},
\; m,m' \in \ibt{1}{M},
\ee
where the infimum  is taken over all signed Borel measures $\mu_1,\ldots,\mu_M$ on $[0,\pi]$.
\elem

\bpf
Let us consider the subspace $\wh{\cU} $ of $\cC[0,\pi]$ defined by $\wh{\cU}  = \{ u \circ \cos, u \in \cU \}$,
for which $(\wh{u}_1,\ldots,\wh{u}_M)$ is a basis.
We readily verify that,
if $P$ is a projection from $\cC[-1,1]$ onto $\cU$,
then 
\be 
Q: g \in \cC[0,\pi] \mapsto P(g \circ \arccos) \circ \cos \in \wh{\cU}
\ee
defines a projection from $\cC[0,\pi]$ onto $\wh{\cU}$ satisfying $\|Q\|_{\infty \to \infty} = \|P\|_{\infty \to \infty}$.
Conversely, we also see that if $Q$ is a projection from $\cC[0,\pi]$ onto $\wh{\cU}$,
then 
\be
P: f \in \cC[-1,1] \mapsto Q(f \circ \cos) \circ \arccos \in \cU
\ee
defines a projection from $\cC[-1,1]$ onto $\cU$
satisfying $\|P\|_{\infty \to \infty} = \|Q\|_{\infty \to \infty}$.
This implies that the projection constant of $\cU$ (relative to $\cC[-1,1]$)
equals the projection constant of $\wh{\cU}$ (relative to $\cC[0,\pi]$),
i.e.,
\be
\label{Ori2}
\la(\cU)  = \la(\wh{\cU}) =
\min_{Q} \; \|Q\|_{\infty \to \infty}
\quad
\mbox{s.to $Q$ being a projection from $\cC[0,\pi]$ onto $\wh{\cU}$}.
\ee
The derivation of \eqref{Transformed} from \eqref{Ori2}
is similar to the derivation of \eqref{Original} from \eqref{Ori1}.
\epf

Lemma \ref{LemTransformed} clearly yields a lower bound for the projection constant if we replace the maximum over the interval $[0,\pi]$ by the maximum over a discretization grid $\theta_1,\ldots,\theta_L$ of $[0,\pi]$,
i.e., we have
\be
\label{Transformed2}
\la(\cU) \ge  \inf_{\mu_1,\ldots,\mu_M}
\max_{\ell \in \lbt 1 : L \rbt} \int_{0}^\pi \left| \sum_{m=1}^M \wh{u}_m(\theta_\ell) d\mu_m \right|
\quad 
\mbox{s.to }
\int_{0}^\pi \wh{u}_{m'} d\mu_{m} = \delta_{m,m'},
\; m,m' \in \ibt{1}{M},
\ee
where the infimum  is taken over all signed Borel measures $\mu_1,\ldots,\mu_M$ on $[0,\pi]$.
Our next step consists in recasting the latter minimization problem 
so as to involve only (nonnegative) Borel measures instead of signed Borel measures.
Although the following observation may seem obvious,
we make an extra effort to verify it fully.

\begin{samepage}
\blem
Given a basis $(u_1,\ldots,u_M)$ for a subspace $\cU$ of $\cC[-1,1]$, 
let $\wh{u}_1,\ldots,\wh{u}_M$ still denote the functions 
$u_1 \circ \cos, \ldots,  u_M \circ \cos$.
If $\theta_1,\ldots,\theta_L \in [0,\pi]$,
then
\begin{align}
\label{InfPBM}
\la(\cU) \ge
\inf_{\substack{\mu_1^\pm,\ldots,\mu_M^\pm\\ \nu_1^\pm, \ldots, \nu_L^\pm}} \max_{\ell \in \lbt 1 : L \rbt} 
\int_{0}^\pi \left( d\nu_\ell^+ + d\nu_\ell^- \right)
& & \mbox{s.to} & 
\int_0^\pi \wh{u}_{m'} (d\mu_m^+ - d \mu_m^-) = \delta_{m,m'}, \; m,m' \in \ibt{1}{M},\\
\nonumber
& & \mbox{} &
\sum_{m = 1}^M \wh{u}_m(\theta_\ell) (\mu_m^+ - \mu_m^-) = \nu_\ell^+ - \nu_\ell^-, \; \ell \in \ibt{1}{L},
\end{align}
where the infimum is taken over all (nonnegative) Borel measures
$\mu_1^\pm,\ldots,\mu_M^\pm, \nu_1^\pm, \ldots, \nu_L^\pm$ on $[0,\pi]$.
\elem
\end{samepage}

\bpf
Let $\alpha$ be the value of the infimum in \eqref{Transformed2} and let $\beta$ be the value of the infimum in \eqref{InfPBM}.
To prove that $\alpha \le \beta$,
we consider miminizers $\mu_1^\pm,\ldots,\mu_M^\pm,\nu_1^\pm,\ldots,\nu_L^\pm$ for the problem \eqref{InfPBM}.
By virtue of the first constraint in \eqref{InfPBM},
 the measures $\mu_m^+ - \mu_m^-$, $m \in \ibt{1}{M}$,
are feasible for the problem \eqref{Transformed2}, 
so that, using the second constraint in \eqref{InfPBM},
\begin{align}
\alpha & \le \max_{\ell \in \lbt 1: L \rbt} \int_0^\pi \left| \sum_{m=1}^M \wh{u}_m(\theta_\ell) (d\mu_m^+ - d \mu_m^-) \right|
= \max_{\ell \in \lbt 1: L \rbt} \int_0^\pi |d\nu_\ell^+ - d\nu_\ell^- | \\
\nonumber
& \le \max_{\ell \in \lbt 1: L \rbt} \int_0^\pi (d\nu_\ell^+ + d\nu_\ell^-)
= \beta.
\end{align}
To prove that $\beta \le \alpha$,
we consider minimizers $\mu_1\ldots,\mu_M$ for the problem \eqref{Transformed2}.
Then we write the Jordan decomposition of each $\mu_m$, $m \in \ibt{1}{M}$,
as $\mu_m = \mu_m^+ - \mu_m^-$ for some (nonnegative) Borel measures $\mu_m^\pm$ satisfying $\mu_m^+ \perp \mu_m^-$.
For each $\ell \in \ibt{1}{L}$, we also write the Jordan decomposition of 
$\nu_\ell := \sum_{m=1}^M \wh{u}_m(\theta_\ell) \mu_m$ as $\nu_\ell= \nu_\ell^+ - \nu_\ell^-$
for some (nonnegative) Borel measures $\nu_\ell^\pm$ satisfying $\nu_\ell^+ \perp \nu_\ell^-$.
In particular, we have $\int_0^\pi |d\nu_\ell| = \int_0^\pi (d\nu_\ell^+ + d\nu_\ell^-)$.
Then, noticing that the Borel measures $\mu_1^\pm,\ldots,\mu_M^\pm,\nu_1^\pm, \ldots, \nu_L^\pm$
are feasible for the problem \eqref{InfPBM}, we obtain
\be
\beta \le \max_{\ell \in \lbt 1: L \rbt} \int_0^\pi (d\nu_\ell^+ + d\nu_\ell^-)
= \max_{\ell \in \lbt 1: L \rbt} \int_0^\pi |d\nu_\ell|
= \max_{\ell \in \lbt 1: L \rbt} \int_0^\pi \left| \sum_{m=1}^M \wh{u}_m(\theta_\ell) d\mu_m \right|
= \alpha.
\ee
The proof is now complete.
\epf

The expression in \eqref{InfPBM} is still not directly exploitable due to the infinite-dimensionality of the Borel measures.
We resolve this issue by substituting these measures by their sequences of moments
and then by truncating the moment constraints.
With fewer constraints,
a smaller value for the minimum is produced.
At the same time, since the moments discarded in the constraints do not occur in the objective function either,
they can be removed altogether
to create a finite-dimensional semidefinite program.
We make all of this precise in the proof of the following result,
which presents  the awaited computable lower bound.
Below, the notation ${\rm Toep}_S(y)$ stands for the $S \times S$ Toeplitz matrix constructed from the first $S$ components of a sequence $(y_k)_{k \ge 0}$,
i.e.,
\be
 {\rm Toep}_S(y) = 
\bbmx 
y_0 & y_1 & \cdots & \cdots & y_{S-1}\\
y_{1} & y_{0} & y_{1} & & \vdots \\
\vdots & y_{1} & \ddots & \ddots & \vdots\\
\vdots & & \ddots & \ddots & y_{1}\\
y_{S-1} & \cdots & \cdots & y_{1} & y_{0}
 \ebmx .
\ee

\bprop
\label{PropFinalLB}
Let $\cU$ be a subspace of $\cC[-1,1]$ consisting of polynomials of degree at most $ d$.
Let   $(u_1,\ldots,u_M)$ be a basis for $\cU$,
whose elements have the Chebyshev expansions
\be
\label{CoefUinT}
u_m = \sum_{k=0}^d U_{k,m} T_k,
\qquad m \in \ibt{1}{M}.
\ee
For an integer $S > d$ and for points $\theta_1,\ldots,\theta_L  \in [0,\pi]$,
one has
\begin{align}
\label{FinalLB}
\la(\cU) \ge \hspace{-2mm}
\inf_{\substack{y_1^\pm,\ldots,y_M^\pm\\ z_1^\pm, \ldots, z_L^\pm}} \max_{\ell \in \lbt 1 : L \rbt} 
(z_{\ell,0}^+ + z_{\ell,0}^- )
& & \mbox{s.to } & 
\sum_{k=0}^d U_{k,m'} (y^+_{m,k} - y^-_{m,k}) = \delta_{m,m'}, \; m,m' \in \ibt{1}{M},\\
\nonumber
& &  &
\sum_{m = 1}^M u_m(\cos(\theta_\ell)) (y_m^+ - y_m^-) = z_\ell^+ - z_\ell^-, \; \ell \in \ibt{1}{L},\\
\nonumber
& & & {\rm Toep}_S(y^\pm_m) \succeq 0, \; m \in \ibt{1}{M},
\\
\nonumber
& & & 
{\rm Toep}_S(z^\pm_\ell) \succeq 0, \; \ell \in \ibt{1}{L},
\end{align}
where the infimum is taken over all vectors $y_1^\pm,\ldots,y_M^\pm,z_1^\pm, \ldots, z_L^\pm \in \bR^S$ indexed from $0$ to $S-1$.
\eprop 

\bpf
In the optimization program \eqref{InfPBM}, we substitute the Borel measures $\mu_1^\pm,\ldots,\mu_M^\pm,\mu_1^\pm,\ldots,\nu_L^\pm$
by their infinite sequences $y_1^\pm,\ldots,y_M^\pm,z_1^\pm,\ldots,z_L^\pm$
of moments,
identified as
\be
\label{TrigMom}
y^\pm_{m,k} = \int_0^\pi \cos(k \theta) d\mu^\pm_m(\theta),
\qquad
z^\pm_{\ell,k} = \int_0^\pi \cos(k \theta) d\nu_\ell^\pm(\theta),
\qquad
k \ge 0,
\ee
to reach an equivalent optimization program featuring the objective function $\max_{\ell \in \lbt 1: L \rbt}(z^+_{\ell,0} + z^-_{\ell,0})$.
As for the constraints (duality, consistency, moments),
the first constraint in \eqref{InfPBM} becomes the first constraint in \eqref{FinalLB} by virtue of
\begin{align}
\int_{0}^\pi \wh{u}_{m'}(\theta) (d\mu_m^+(\theta) - d\mu_m^-(\theta))
& =  \int_0^\pi \sum_{k=0}^d U_{k,m'} T_k(\cos(\theta)) (d\mu_m^+(\theta) - d\mu_m^-(\theta))\\
\nonumber
& =  \sum_{k=0}^d U_{k,m'} \int_0^\pi \cos(k \theta) (d\mu_m^+(\theta) - d\mu_m^-(\theta))\\
\nonumber
& =   \sum_{k=0}^d U_{k,m'} (y^+_{m,k} - y^-_{m,k});
\end{align}
the second constraint in \eqref{InfPBM} is clearly equivalent to the second constraint in \eqref{FinalLB};
while the fact that we are dealing with sequences of moments is reflected by the semidefinite conditions ${\rm Toep}_\infty(y_m^\pm) \succeq 0$, $m \in \ibt{1}{M}$, and ${\rm Toep}_\infty(z_\ell^\pm) \succeq 0$, $\ell \in \ibt{1}{L}$.
We now relax the last two sets of constraints by simply imposing,
for $\ell \in \ibt{1}{L}$,
\be
\sum_{m = 1}^M u_m(\cos(\theta_\ell)) (y_{m,k}^+ - y_{m,k}^-) = z_{\ell,k}^+ - z_{\ell,k}^-,
\qquad
k \in \ibt{0}{S \hspace{-1mm} - \hspace{-1mm}1},
\ee
as well as, for $m \in \ibt{1}{M}$ and $\ell \in \ibt{1}{L}$,
\be
{\rm Toep}_S(y_m^\pm) 
 \succeq 0
\qquad \mbox{and} \qquad
{\rm Toep}_S(z_\ell^\pm) \succeq 0.
\ee
This leads to a smaller minimum value for the optimization program.
And since the relaxed program only involves moments up to order $S-1$,
we can restrict the minimization to the finite sequences 
$(y_{m,k}^\pm)_{k=0}^{S-1}$, $m \in \ibt{1}{M}$,
and $(z_{\ell,k}^\pm)_{k=0}^{S-1}$, $\ell \in \ibt{1}{L}$,
hence yielding the computable lower bound stated in \eqref{FinalLB}.
\epf

We close this section by highlighting how the lower bound from Proposition \ref{PropFinalLB}
is expressed as a semidefinite program.
For this purpose,
besides the matrix $U \in \bR^{(d+1) \times M}$
containing the coefficients of $u_1,\ldots,u_M$ in the Chebyshev system $(T_0,\ldots,T_d)$, as indicated in \eqref{CoefUinT},
we also consider the collocation matrix $W \in \bR^{L \times M}$ with entries
\be
\label{DefW2}
W_{\ell,m} = u_m(\cos(\theta_\ell)),
\qquad
\ell \in \ibt{1}{L},
\; m \in \ibt{1}{M}.
\ee
By further introducing matrices $Y^{\pm} \in \bR^{M \times S}$
with rows $y_1^\pm,\ldots,y_M^\pm \in \bR^S$
and $Z^{\pm} \in \bR^{L \times S}$
with rows $z_1^\pm,\ldots,z_L^\pm \in \bR^S$,
as well as a slack variable $c \in \bR$ 
(such that $z^+_{\ell,0} + z^-_{\ell,0} \le c$ for all $\ell$),
the optimization program in \eqref{FinalLB} 
takes the easily implementable form below (where some convenient {\sc matlab} notation is used).\footnote{The grid points $\cos(\theta_1),\ldots,\cos(\theta_L)$ could be added as inputs --- by default, we chose them to be Chebyshev zeros.\label{FNgp}}

\ovalbox{
\begin{minipage}{0.97\textwidth}
\medskip
\begin{center}
\textbf{Lower bound for the projection constant}\vspace{-2mm}\\
%\rule{4.6in}{.8pt}
\rule{0.96\textwidth}{.8pt}
\end{center} 
Inputs: basis for a space $\cU$ of polynomials of degree $\le d$, parameters $S>d$ and $L$. 
\be
\label{ProgLB}
\la(\cU) 
\ge 
\min_{\substack{Y^\pm \in \bR^{M \times S}\\ Z^\pm \in \bR^{L \times S}\\ c \in \bR }}
\; \; c
 \qquad \mbox{s.to \;\;}
\left\{ 
\bmx
(Y^+(:,1\hspace{-1mm}: \hspace{-1mm}d+1)-Y^-(:,1\hspace{-1mm}: \hspace{-1mm}d+1)) \, U= I_M,  \\
 W \, (Y^+ - Y^-)  = Z^+ - Z^-, \hfill \\
{\rm Toep}_S(Y^\pm(m,:)) \succeq 0,
\quad  m \in \ibt{1}{M}, \hfill\\
{\rm Toep}_S(Z^\pm(\ell,:)) \succeq 0,
\quad  \ell \in \ibt{1}{L}, \hfill \\
Z^+(:,1) + Z^-(:,1) \le c , \hfill
\emx
\right.
\ee
where the matrices $U$ and $V$ are defined in \eqref{CoefUinT} and \eqref{DefW2}.
\medskip
\end{minipage}
}

\section{Exploiting the symmetry of minimal projections}
\label{SecSym}

The programs highlighted in \eqref{ProgUB} and \eqref{ProgLB} are computationally demanding for large values of the parameters $K$, $L$, and $S$,
so any property that can reduce their complexity should be exploited.
We shall capitalize on a certain symmetry of minimal projections.
To this end, we assume from now on that the subspace $\cU$ of $\cC[-1,1]$ satisfies
\be
\label{USym}
u \in \cU \imp u(-\cdot) \in \cU,
\ee
where $u(-\cdot)$ evidently denotes the function $x \in [-1,1] \mapsto u(-x) \in \bR$.
Under this assumption, it is known that there exists a minimal projection $P$ from $\cC[-1,1]$ onto $\cU$ which is symmetric, in the sense that
\be
\label{SymCond}
P(f(-\cdot))= (P(f))(-\cdot)
\qquad \mbox{for all $f \in \cC[-1,1]$}.
\ee
This fact has the following implication whenever assumption \eqref{USym} holds.\footnote{Under assumption \eqref{USym}, one can verify that there exists a basis $(u_1,\ldots,u_M)$ whose elements are either even or odd functions
(verify, for instance, that $\cU = \cU^{\rm e} \overset{\perp}{\oplus} \cU^{\rm o}$,
where $\cU^{\rm e/o} := \{ u \in \cU: u \mbox{ is an even/odd function} \}$,
and concatenate a basis for $\cU^{\rm e}$ with a basis for $\cU^{\rm o}$).}

\bprop
Let $\cU$ be a subspace of $\cC[-1,1]$ and let $(u_1^{\rm e},\ldots,u_{M_{\rm e}}^{\rm e},u_1^{\rm o},\ldots,u_{M_{\rm o}}^{\rm o})$ be a basis for~$\cU$ arranged in such a way that the $u_m^{\rm e}$ are even functions and the $u_m^{\rm o}$ are odd functions.
Considering functions $\wt{u}_1^{\rm e},\ldots,\wt{u}_{M_{\rm e}}^{\rm e},\wt{u}_1^{\rm o},\ldots,\wt{u}_{M_{\rm o}}^{\rm o}$ defined on $[-1,1]$ by
\be
\wt{u}^{\rm e/o}_m(t) := u^{\rm e/o}_m \left( \f{t+1}{2} \right),
\qquad t \in [-1,1],
\ee
the projection constant of $\cU$ can be expressed as
\begin{align}
\label{NewExpr}
\la(\cU) = & \inf_{\substack{\wt{\mu}_1^{\rm e},\ldots,\wt{\mu}_{M_{\rm e}}^{\rm e}\\ \wt{\mu}_1^{\rm o},\ldots,\wt{\mu}_{M_{\rm o}}^{\rm o}}}
 \max_{x \in [0,1]}
 \int_{-1}^1
\max  \left\{  
\left| \sum_{m=1}^{M_{\rm e}} u_m^{\rm e}( x ) d\wt{\mu}_m^{\rm e} \right|,
\left| \sum_{m=1}^{M_{\rm o}} u_m^{\rm o}( x ) d\wt{\mu}_m^{\rm o} \right|
\right\}\\
\nonumber
\mbox{s.to} & \int_{-1}^1 \wt{u}_{m'}^{\rm e} d\wt{\mu}_m^{\rm e} = \delta_{m,m'}, \; m,m' \in \ibt{1}{M_{\rm e}},
\qquad
 \int_{-1}^1 \wt{u}_{m'}^{\rm o} d\wt{\mu}_m^{\rm o} = \delta_{m,m'}, \; m,m' \in \ibt{1}{M_{\rm o}},
\end{align}
where the infimum is taken over all signed Borel measures $\wt{\mu}_1^{\rm e},\ldots,\wt{\mu}_{M_{\rm e}}^{\rm e},\wt{\mu}_1^{\rm o},\ldots,\wt{\mu}_{M_{\rm o}}^{\rm o}$ on $[-1,1]$.
\eprop

\bpf
Let us consider a symmetric projection from $\cC[-1,1]$ onto $\cU$ written as
\be
P(f) = \sum_{m=1}^{M_{\rm e}} \eta^{\rm e}_m(f) u_m^{\rm e}
+ \sum_{m=1}^{M_{\rm o}} \eta^{\rm o}_m(f) u_m^{\rm o},
\qquad f \in \cC[-1,1].
\ee
The condition \eqref{SymCond} is readily seen to be equivalent to the conditions
\be
\label{SymEta}
\eta_m^{\rm e}(f(-\cdot)) =  \eta_m^{\rm e}(f),
\qquad
\eta_m^{\rm o}(f(-\cdot)) =  -\eta_m^{\rm o}(f),
\qquad f \in \cC[-1,1],
\ee
which in turn are equivalent, in terms of measures $\mu_1^{\rm e}, \ldots, \mu_{M_{\rm e}}^{\rm e},\mu_1^{\rm o},\ldots,\mu_{M_{\rm o}}^{\rm o}$
representing the linear functionals $\eta_1^{\rm e}, \ldots, \eta_{M_{\rm e}}^{\rm e},\eta_1^{\rm o},\ldots,\eta_{M_{\rm o}}^{\rm o}$, to the conditions
\be
\label{SymMu}
d\mu_m^{\rm e}(- \cdot) = d\mu_m^{\rm e},
\qquad
d\mu_m^{\rm o}(- \cdot) = - d\mu_m^{\rm o}.
\ee
Then, the norm of the projection $P$ satisfies
\begin{align}
\|P\|_{\infty \to \infty}
& = \max_{x \in [-1,1]} \int_{-1}^1 \left| \sum_{m=1}^{M_{\rm e}} u_m^{\rm e}(x) d\mu_m^{\rm e} + \sum_{m=1}^{M_{\rm o}} u_m^{\rm o}(x) d\mu_m^{\rm o} \right|\\
\nonumber & = \max_{x \in [-1,1]} \int_{-1}^1 \left| \sum_{m=1}^{M_{\rm e}} u_m^{\rm e}(x) d\mu_m^{\rm e} - \sum_{m=1}^{M_{\rm o}} u_m^{\rm o}(x) d\mu_m^{\rm o} \right|\\
\nonumber & = \max_{x \in [-1,1]} \int_{-1}^1 \max \left\{ \left| \sum_{m=1}^{M_{\rm e}} u_m^{\rm e}(x) d\mu_m^{\rm e}\right| , \left| \sum_{m=1}^{M_{\rm o}} u_m^{\rm o}(x) d\mu_m^{\rm o} \right| \right\}, 
\end{align}
where we have used the identity $(|a+b|+|a-b|) / 2 = \max\{ |a|, |b| \}$.
Noticing the invariance of the above expression under the change $x \leftrightarrow -x$
and taking \eqref{SymMu} into account,
we can further write
\begin{align} 
\label{ExpNormSym}
\|P\|_{\infty \to \infty}
& = \max_{x \in [0,1]} 2 \int_{0}^1 \max \left\{ \left| \sum_{m=1}^{M_{\rm e}} u_m^{\rm e}(x) d\mu_m^{\rm e}\right| , \left| \sum_{m=1}^{M_{\rm o}} u_m^{\rm o}(x) d\mu_m^{\rm o} \right| \right\}\\
\nonumber
& = \max_{x \in [0,1]} \phantom{2} \int_{-1}^1 \max \left\{ \left| \sum_{m=1}^{M_{\rm e}} u_m^{\rm e}(x) d\wt{\mu}_m^{\rm e}\right| , \left| \sum_{m=1}^{M_{\rm o}} u_m^{\rm o}(x) d\wt{\mu}_m^{\rm o} \right| \right\}, 
\end{align}
where the measures $\wt{\mu}_m^{\rm e/o}$ simply represent the restrictions to $[0,1]$ of the measures $\mu_m^{\rm e/o}$ that have been transposed to $[-1,1]$,
i.e., they are obtained through the identification
\be
d\wt{\mu}_m^{\rm e/o}(t) = 2 d \mu_m^{\rm e/o}(\tau),
\qquad  t \in [-1,1], \tau \in [0,1]
\mbox{ being linked via }
t = 2 \tau - 1, \,
\tau = \f{t+1}{2}. 
\ee
Thanks to \eqref{SymMu},
we can take the $\wt{\mu}_m^{\rm e/o}$ as new optimization variables in the minimization of the norm of a symmetric projection,
whose expression \eqref{ExpNormSym} is the objective function in \eqref{NewExpr}.
We now just have to impose the appropriate duality constraints making $P$ a projection onto~$\cU$.
Among them,
the constraints $\eta_m^{\rm e/o}(u_{m'}^{\rm o/e}) = 0$ are automatically fulfilled,
while the constraints $\eta_m^{\rm e/o}(u_{m'}^{\rm e/o}) = \delta_{m,m'}$ reduce to
\be
\delta_{m,m'} = \int_{-1}^1 u_{m'}^{\rm e/o} d\mu_{m}^{\rm e/o}
= 2 \int_{0}^1 u_{m'}^{\rm e/o}(\tau) d\mu_{m}^{\rm e/o}(\tau)
= \int_{-1}^1 \wt{u}_{m'}^{\rm e/o} (t) d\wt{\mu}_m^{\rm e/o}(t).
\ee
These are indeed the constraints in \eqref{NewExpr},
so the proposition is proved.
\epf

We proceed by highlighting the computable bounds on the projection constant generated by the reformulation \eqref{NewExpr}. 
Much of the ingredients for deriving these bounds are similar to the ones presented in Sections \ref{SecUB} and \ref{SecLB},
so we do expand on details at all.

\subsection{Implication for the upper bound}

As in Section \ref{SecUB}, we first derive an upper bound for the projection constant by minimizing only over linear combinations 
\be
\wt{\mu}_m^{\rm e/o} = \sum_{k=1}^K A^{\rm e/o}_{m,k} \delta_{v_k}.
\ee
of Dirac measures at $v_1,\ldots,v_K \in [-1,1]$.
Then we again discretize by replacing the maximum over $[0,1]$ by the maximum over
a grid $w_1^+ > \cdots > w_L^+$ consisting of positive Chebyshev zeros 
$w_\ell^+ = \cos(\pi (\ell-1/2)/(2L))$.
Thus,
we arrive at the computable upper bound
\begin{align}
\la(\cU)
\le \rho \times \hspace{-1mm}
\min_{A^{\rm e/o} \in \bR^{M_{\rm e/o} \times K}} 
\max_{\ell \in \lbt 1 : L \rbt } 
\sum_{k=1}^K 
\max & \left\{
\left| \sum_{m=1}^{M_{\rm e}} A^{\rm e}_{m,k} u^{\rm e}_m(w_\ell^+) \right|,
\left| \sum_{m=1}^{M_{\rm o}} A^{\rm o}_{m,k} u^{\rm o}_m(w_\ell^+) \right|
\right\}\\
\nonumber
\mbox{s.to } & 
 \sum_{k=1}^K A^{\rm e/o}_{m,k} \wt{u}^{\rm e/o}_{m'}(v_k) = \delta_{m,m'},
 \,  m,m' \in \ibt{1}{M_{\rm e/o}},
\end{align}
where $\rho = \cos \left( (\pi d)/(4 L ) \right)^{-1}$.
To transform the latter into a linear program,
we introduce the collocation matrices $V^{\rm e/o} \in \bR^{K \times M_{\rm e/o}}$
and $W^{\rm e/o} \in \bR^{L \times M_{\rm e/o}}$ defined by
\be
\label{DefVW3}
V^{\rm e/o}_{k,m} = \wt{u}^{\rm e/o}_m(v_k)
\qquad \mbox{and} \qquad
W^{\rm e/o}_{\ell,m} = u^{\rm e/o}_m(w_{\ell}^+),
\qquad
k \in \ibt{1}{K}, \;
\ell \in \ibt{1}{L}, \;
m \in \ibt{1}{M_{\rm e/o}},
\ee
as well as slack variables through $B \in \bR^{L \times K}$ and $c \in \bR$.
All in all, we obtain the following implementable form of the upper bound.

\ovalbox{
\begin{minipage}{0.97\textwidth}
\medskip
\begin{center}
\textbf{Upper bound for the projection constant --- symmetry exploited}\vspace{-2mm}\\
%\rule{4.6in}{.8pt}
\rule{0.96\textwidth}{.8pt}
\end{center} 
Inputs: basis for a symmetric space $\cU$ of polynomials of degree $\le d$, parameters $S>d$ and $L$. 
\be
\label{ProgUBSym}
\la(\cU) 
\le \cos \left( \df{\pi}{4}\df{d}{L} \right)^{-1}
\times
\min_{\substack{ A^{\rm e/o} \in \bR^{M_{\rm e/o} \times K}
 \\
B \in \bR^{L \times K} \\ c \in \bR}} 
\;  c 
\qquad
 \mbox{s.to \; \;} 
 \left\{ 
 \bmx
A^{\rm e/o} V^{\rm e/o} = I_{M_{\rm e/o}},  \hfill \\
-B \le W^{\rm e/o} A^{\rm e/o} \le B, \hfill \\
B {\bf 1} \le c {\bf 1}, \hfill
\emx
\right.
\ee
where the matrices $V^{\rm e/o}$ and $W^{\rm e/o}$ are defined in \eqref{DefVW3}.
\medskip
\end{minipage}
} 

\subsection{Implication for the lower bound}

As in Section \ref{SecLB}, the minimization program \eqref{NewExpr} is first transformed to make it amenable to the trigonometric moment problem.
The reformulation will involve a maximum over $[0,\pi/2]$,
which is lower bounded by the maximum over a grid $\theta_1^+,\ldots,\theta_L^+$.
With $w_\ell^+ := \cos(\theta_\ell^+) \in [0,1]$, we obtain
\begin{align}
\la(\cU) \ge \inf_{\substack{\wt{\mu}_1^{\rm e},\ldots,\wt{\mu}_{M_{\rm e}}^{\rm e}\\ \wt{\mu}_1^{\rm o},\ldots,\wt{\mu}_{M_{\rm o}}^{\rm o}}}
\max_{\ell \in \lbt 1: L \rbt }
 \int_{0}^\pi
\max & \left\{  
\left| \sum_{m=1}^{M_{\rm e}} u_m^{\rm e}(w_\ell^+) d\wt{\mu}_m^{\rm e} \right|,
\left| \sum_{m=1}^{M_{\rm o}} u_m^{\rm o}(w_\ell^+) d\wt{\mu}_m^{\rm o} \right|
\right\}\\
\nonumber
\mbox{s.to} & \int_{0}^\pi \wt{u}_{m'}^{\rm e/o} (\cos( \theta) ) d\wt{\mu}_m^{\rm e/o}(\theta) = \delta_{m,m'}, \quad m,m' \in \ibt{1}{M_{\rm e/o}},
\end{align}
where the infimum is taken over all signed Borel measures $\wt{\mu}_1^{\rm e},\ldots,\wt{\mu}_{M_{\rm e}}^{\rm e},\wt{\mu}_1^{\rm o},\ldots,\wt{\mu}_{M_{\rm o}}^{\rm o}$ on $[0,\pi]$.
\mbox{Deviating} slightly from the earlier strategy,
we now introduce as slack variables some \mbox{(nonnegative)} Borel measures $\nu_1,\ldots,\nu_L$
 satisfying
\be 
\nu_\ell \ge 
\max \left\{  
\left| \sum_{m=1}^{M_{\rm e}} u_m^{\rm e}(w_\ell^+) \wt{\mu}_m^{\rm e} \right|,
\left| \sum_{m=1}^{M_{\rm o}} u_m^{\rm o}(w_\ell^+) \wt{\mu}_m^{\rm o} \right|
\right\},
\qquad \mbox{i.e.,} \qquad
\nu_\ell \ge \pm 
\sum_{m=1}^{M_{\rm e/o}} u_m^{\rm e/o}(w_\ell^+) \wt{\mu}_m^{\rm e/o}.
\ee 
Writing the Chebyshev expansions of $\wt{u}_1^{\rm e},\ldots, \wt{u}_{M_{\rm e}}^{\rm e}, \wt{u}_1^{\rm o},\ldots, \wt{u}_{M_{\rm o}}^{\rm o}$
as
\be
\label{DefU4}
\wt{u}_m^{\rm e/o} = \sum_{k=0}^d \wt{U}^{\rm e/o}_{k,m} T_k,
\qquad m \in \ibt{1}{M_{\rm e/o}},
\ee
and substituting the measures $\wt{\mu}_m^{\rm e/o}$ and $\nu_\ell$
by their sequences $y_m^{\rm e/o}$ and $z_\ell$ of moments 
yields to
\begin{align}
\la(\cU) \ge \hspace{-1mm}
\inf_{\substack{y_1^{\rm e/o},\ldots,y_{M_{\rm e/o}}^{\rm e/o}\\ z_1,\ldots, z_L}}
\max_{\ell \in \lbt 1: L \rbt } z_{\ell,0}
\quad \mbox{s.to \;  } & 
\sum_{k=0}^{d} \wt{U}^{\rm e/o}_{k,m'} y_{m,k}^{\rm e/o} = \delta_{m,m'},
\; m,m' \in \ibt{1}{M_{\rm e/o}},\\
\nonumber
& {\rm Toep}_\infty(z_\ell) \succeq \pm {\rm Toep}_\infty \left( \sum_{m=1}^{M_{\rm e/o}} u_m^{\rm e/o} ( w_\ell^+ ) y_m^{\rm e/o} \right),
\; \ell \in \ibt{1}{L},
\end{align}
where the infinum is taken over all infinite sequences $y_1^{\rm e/o},\ldots,y_{M_{\rm e/o}}^{\rm e/o}, z_1,\ldots, z_L \in \bR^{\bN}$.
The infinite semidefinite constraints are now truncated to a level $S>d$,
producing a lower bound involving only the finite sequences of moments $(y_{m,k}^{\rm e/o})_{k=0}^{S-1}$, $m \in \ibt{1}{M_{\rm e/o}}$, and $(z_{\ell,k})_{k=0}^{S-1}$, $\ell \in \ibt{1}{L}$.
Finally, in order to state the corresponding semidefinite program in an easily implementable form,
we gather these moments in matrices $Y^{\rm{e/o}} \in \bR^{M_{\rm e/o}\times S}$ and $Z \in \bR^{L \times S}$,
and we define collocation matrices $W^{\rm e/o} \in \bR^{L \times M_{\rm e/o}}$ with entries
\be
\label{DefW4}
W^{\rm e/o}_{\ell,m} = u^{\rm e/o}_m(w_\ell^+),
\qquad
\ell \in \ibt{1}{L},
\; m \in \ibt{1}{M_{\rm e/o}}.
\ee
After introducing one last slack variable $c \in \bR$,
we arrive at the following form of the computable lower bound
(where some convenient {\sc matlab} notation is again used).\footnote{In a similar spirit to footnote \ref{FNgp}, 
the grid points $w_1^+ = \cos(\theta_1), \ldots,w_L^+ =\cos(\theta_L)$ could be added as inputs --- by default, we chose them to be positive Chebyshev zeros.}

\ovalbox{
\begin{minipage}{0.97\textwidth}
\medskip
\begin{center}
\textbf{Lower bound for the projection constant --- symmetry exploited}\vspace{-2mm}\\
%\rule{4.6in}{.8pt}
\rule{0.96\textwidth}{.8pt}
\end{center} 
Inputs: basis for a symmetric space $\cU$ of polynomials of degree $\le d$, parameters $S>d$ and $L$. 
\be
\label{ProgLBSym}
\la(\cU) 
\ge 
\min_{\substack{Y^{\rm e/o} \in \bR^{M_{\rm e/o} \times S}\\ Z \in \bR^{L \times S}\\ c \in \bR }}
\; c
 \qquad \mbox{s.to \;\;}
\left\{ 
\bmx
Y^{\rm e/o}(:,1\hspace{-1mm}: \hspace{-1mm}d+1) \, \wt{U}^{\rm e/o}= I_{M_{\rm e/o}}, \hfill \\
{\rm Toep}_S(Z(\ell,:)) \succeq 
\pm {\rm Toep}_S(W^{\rm e/o}(\ell,:)Y^{\rm e/o}),
\;   \ell \in \ibt{1}{L}, \hfill \\
Z(:,1) \le c , \hfill
\emx
\right.
\ee
where the matrices $\wt{U}^{\rm e/o}$ and $W^{\rm e/o}$ are defined in \eqref{DefU4} and \eqref{DefW4}.
\medskip
\end{minipage}
}

\section{Computational results}
\label{SecCompRes}

This section gives an account of the experiments carried out using our method for specific spaces of univariate polynomials.
The experiments can be reproduced by downloading the {\sc matlab} file tied to this paper, available on the authors' webpages.
Note that the code relies on {\sf CVX} \cite{cvx},
a {\sc matlab} package for specifying and solving convex programs,
and on {\sf Chebfun}~\cite{Chebfun} for its convenience to deal with Chebyshev expansions.

\subsection{Validation of the code}

\begin{figure}[h]
\center
\subfigure{
\includegraphics[width=0.31\textwidth]{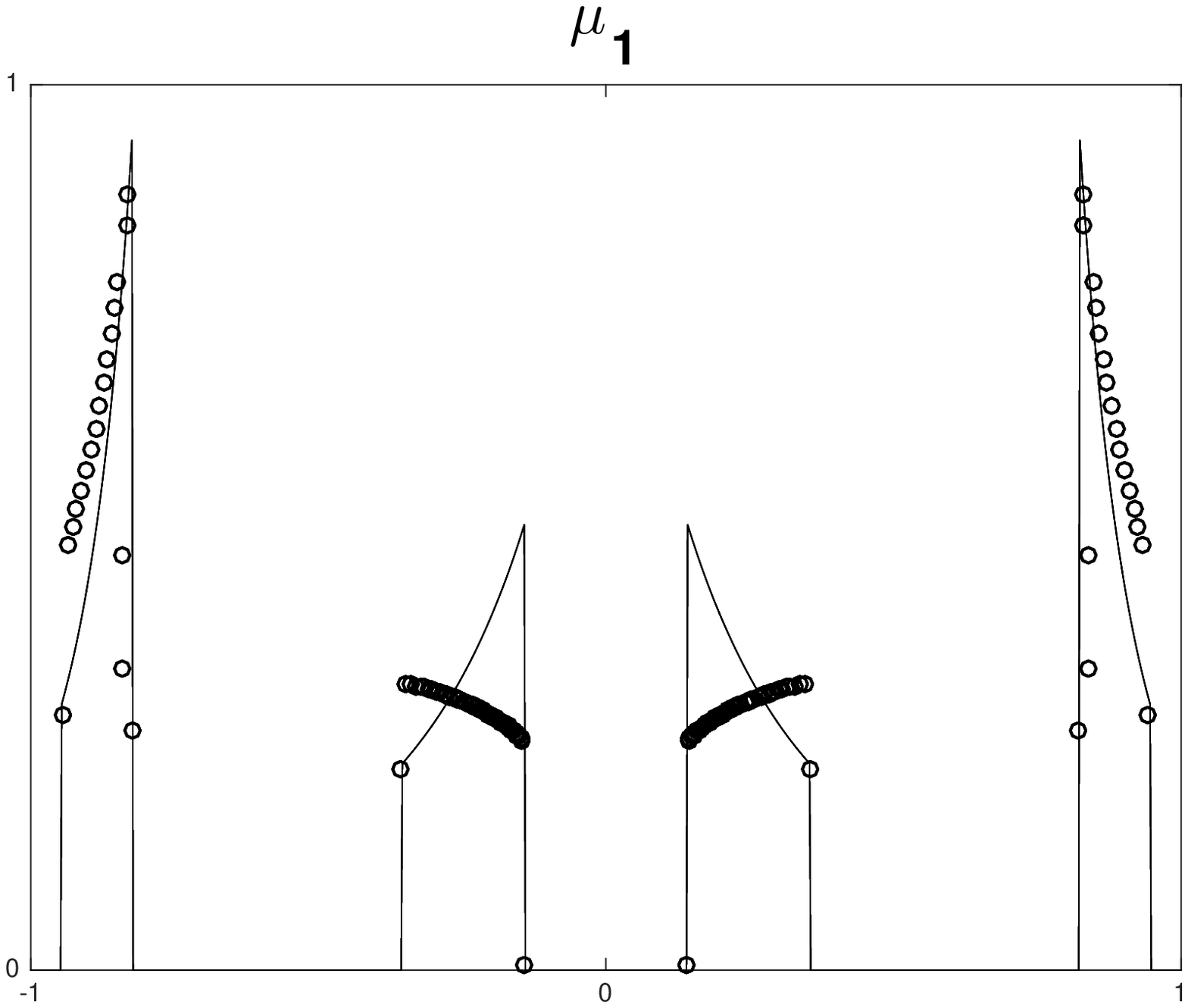}
}
\subfigure{
\includegraphics[width=0.31\textwidth]{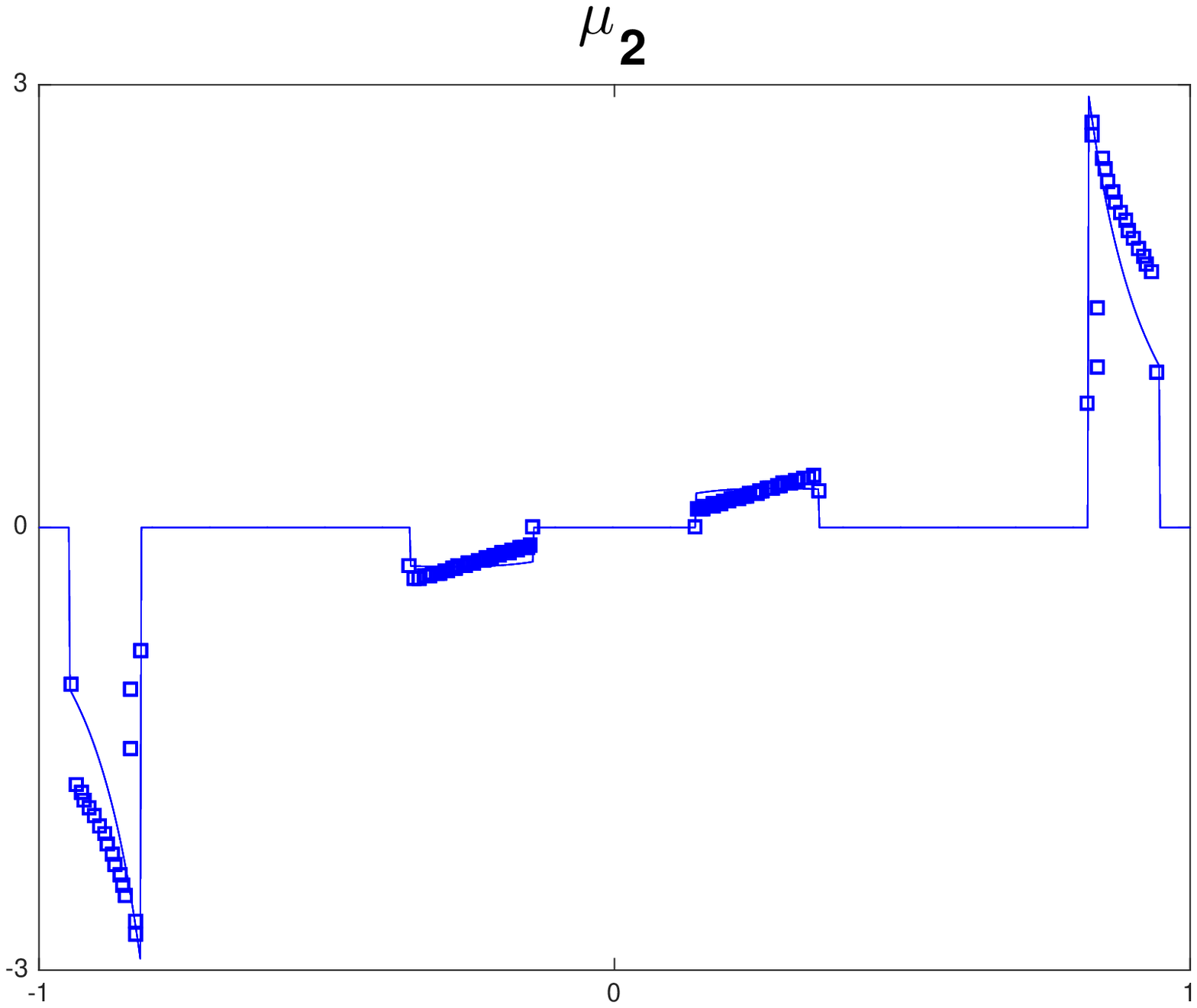}
}
\subfigure{
\includegraphics[width=0.31\textwidth]{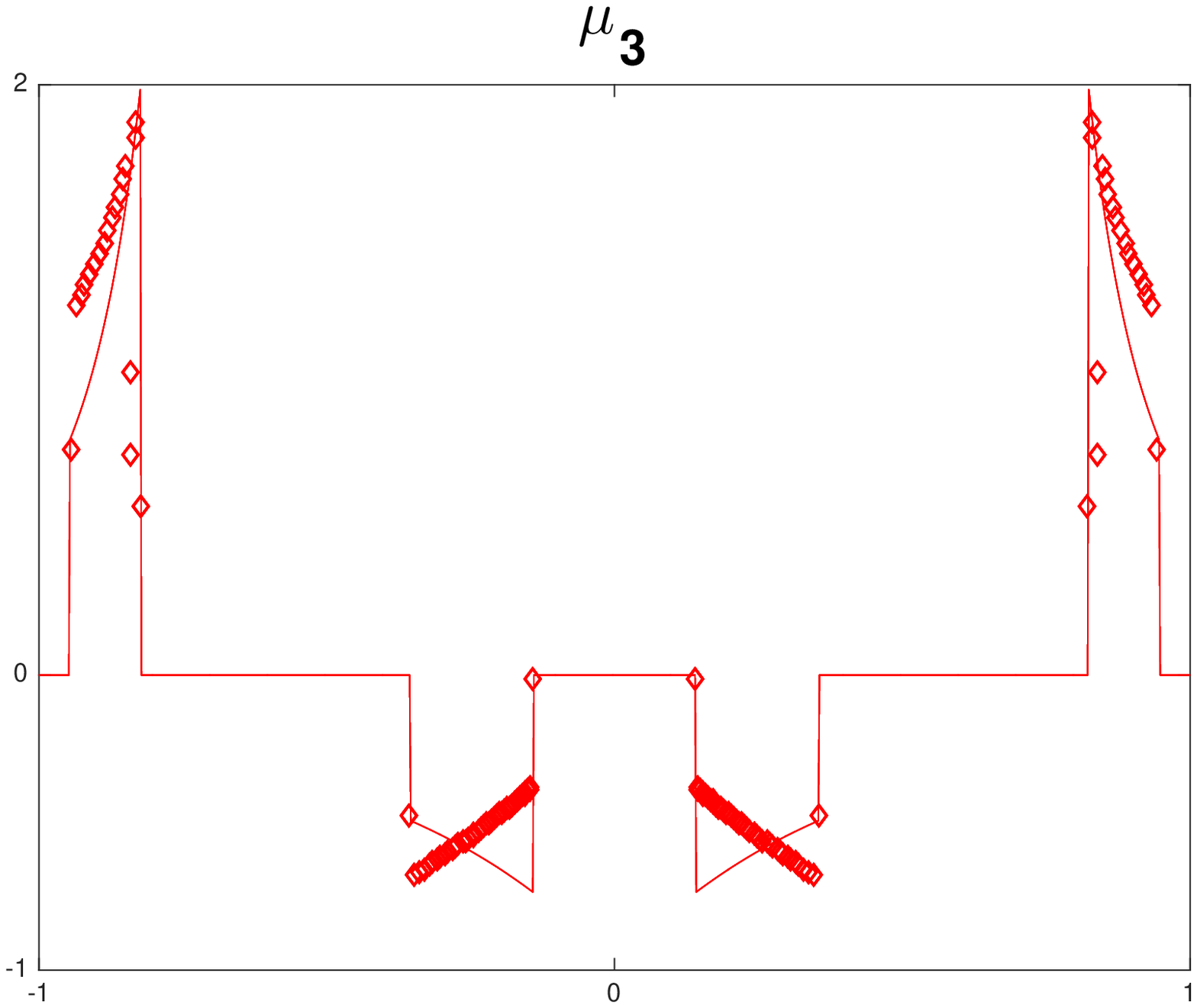}
}
\caption{The nonatomic parts of measures $\mu_1,\mu_2,\mu_3$ associated with the linear functionals $\eta_1,\eta_2,\eta_3$ found in \cite{ChaMet} (continuous lines),
together with their discrete approximations obtained as solutions of \eqref{ProgUBSym},
minus the parts at $-1$, $0$, and $1$ (circles, squares, diamonds).}
\label{Fig1}
\end{figure}

In order to certify the correct implementation of the codes computing the upper and lower bounds \eqref{ProgUBSym} and \eqref{ProgLBSym},
we take as a benchmark the inevitable result \cite{ChaMet} of Chalmers and Metcalf,
who managed to determine analytically the projection constant of the space of quadratic polynomials.
They obtained
\be
\label{ValQuad}
\la(\cP_2) \approx 1.220173064217988\ldots
\ee
and exhibited a minimal projection given by $P(f)(x) = \sum_{m=1}^3 \eta_m(f) x^{m-1}$, 
where
\begin{eqnarray}
\label{CMMeas1}\eta_1(f) = &  A f(-1) + B f(0) + A f(1) & + \int_{I_1} \f{a_1 |t| + b_1}{(1+w_1 |t|)^3} f(t)dt + \int_{I_2} \f{a_2 |t| + b_2}{(1+w_2 |t|)^3} f(t)dt  ,\\
\label{CMMeas2}\eta_2(f) = &  -C f(-1) + C f(1) & + \int_{I_1} \f{c_1t}{(1+w_1 |t|)^3} f(t) dt + \int_{I_2} \f{c_2t}{(1+w_2 |t|)^3} f(t) dt,\\
\label{CMMeas3}\eta_3(f) = & D f(-1) - B f(0) + D f(1) & + \int_{I_1} \f{d_1 |t| - b_1}{(1+w_1 |t|)^3} f(t)dt + \int_{I_2} \f{d_2 |t| - b_2}{(1+w_2 |t|)^3} f(t)dt,
\end{eqnarray}
with $I_1 = [-s_{1,2},-s_{1,1}] \cup [s_{1,1},s_{1,2}]$, $I_2=[-s_{2,2},-s_{2,1}] \cup [s_{2,1} ,s_{2,2}]$, and with parameters $A,B,C,D$, $a_1,b_1,c_1,d_1$, $a_2,b_2,c_2,d_2$, $w_1,w_2$, $s_{1,1}, s_{1,2}, s_{2,1}, s_{2,2}$ determined in the body of \cite{ChaMet}.  
Our code does allow us to retrieve the value \eqref{ValQuad}
up to five digits of accuracy.\footnote{in about five minutes for the upper bound and fifteen minutes for the lower bound,
with the capabilities offered by a laptop computer at the time this article was written.}
Incidentally, our experiment suggests that minimal projections onto the quadratics are not unique,
as illustrated by Figure \ref{Fig1} which superimposes the measures associated to \eqref{CMMeas1}-\eqref{CMMeas2}-\eqref{CMMeas3} and the discretized measures obtained by solving \eqref{ProgUBSym} ---
we have removed the atomic parts at $-1,0,1$,
which were similar.
It is worth noticing, nonetheless, that the supports of all the measures seem to be the same.

Another way validate of our code is offered by the three-dimensional space ${\rm span}\{1,x^2,x^3\}$, i.e., 
the subspace $\cC[-1,1]$ spanned by $ x \mapsto 1$, $x \mapsto x^2$, and $x \mapsto x^3$.
%incorrectly but unambiguously 
Indeed, its projection constant is exactly equal to one.
More generally, 
for any even integer $a>0$ 
and  any odd integer $b>0$,
the space ${\rm span}\{ 1, x^a,x^{a+b} \} = {\rm span}\{1-x^a, x^a(1-x^b),x^a(1+x^b)\}$
admits a projection of norm one.
To see this, we notice that  the interpolating projection at the points $-1$, $0$, and $1$,
which is given by 
\be
P(f)(x) = \f{1}{2} f(-1) x^{a}(1-x^{b}) + f(0) (1-x^{a}) + \f{1}{2} f(1) x^{a}(1+x^{b}),
\ee
satisfies, for all $f \in \cC[-1,1]$ and $x \in [-1,1]$,
\be
|P(f)(x)| \le \max\{ |f(-1)|,|f(0)|,|f(1)| \} \left( \f{1}{2}x^a(1-x^b) + 1-x^a + \f{1}{2} x^a(1+x^b) \right)
\le \|f\|_\infty .
\ee
Our code confirms the value $\la({\rm span}\{1,x^2,x^3\}) = 1$
(not with perfect accuracy, though, which is why we refrain from supplying numerical values with more than five digits in general).

\subsection{Other-three dimensional polynomial spaces}

The code distributed with the {\sc matlab} reproducible can be effortlessly applied to any univariate polynomial space,
so long as it can be executed with parameters large enough for the upper and lower bounds to match up to the desired accuracy.
Five digits of accuracy can typically be achieved for three-dimensional spaces of moderate degree.
The following table summarizes the values of projection constants obtained for several different spaces.\footnote{For the space ${\rm span}\{T_1,T_2,T_3\}$, it was more effective to compute the lower bound using a `symmetrization' of \eqref{ProgLB} more direct than \eqref{ProgLBSym}.
Its implementation is also included in the reproducible.}

\begin{table}[h]
\begin{center}
\begin{tabular}{|c|c|c|c|c|c|c|}
\hline
$\cV$ spanned by  & 
$\{ 1,x,x^3 \}$ & $\{T_0,T_2,T_3\}$ & $\{U_0,U_2,U_3\}$ & $\{ x,x^2,x^3 \}$ & $\{T_1,T_2,T_3\}$ & $\{U_1,U_2,U_3 \}$\\
\hline
$\la(\cV)$ & $\approx 1.4723$ & $\approx 1.4460$ & $\approx 1.1522$ & $\approx 1.3325$ & $\approx 1.4065$ &   $\approx 1.2354$ \\
 \hline
\end{tabular}
\end{center}
\caption{Projection constants of three-dimensional spaces spanned by monomials and Chebyshev polynomials of the first and second kind.}
\end{table}

\subsection{Polynomial spaces of higher dimensions}

Looking at \eqref{ProgUBSym} and \eqref{ProgLBSym},
we see that the dimension $M=M_{\rm e} + M_{\rm o}$
is far from being a factor influencing the computational cost of the optimization programs,
so our code can easily be executed for polynomial spaces of dimension higher than three.
For instance,
we can deal with the spaces of cubic, quartic, and quintic polynomials
and compute their projection constants
with four digits of accuracy as
%\begin{align} 
%\la(\cP_3) &\approx 1.365,\\
%\la(\cP_4) & \approx 1.459,\\
%\la(\cP_5) & \approx 1.538.
%\end{align}
\be 
\la(\cP_3) \approx 1.365,\qquad
\la(\cP_4) \approx 1.459,\qquad
\la(\cP_5) \approx 1.538.
\ee

\begin{samepage}
Obtaining the same accuracy necessitates larger parameters $K$, $L$, $S$
when the dimension increases.
For degree $d>5$, with our modest computational investment,
we could locate the projection constants of the spaces $\cP_d$ of polynomials of degree at most $d$ in the ranges presented in Table \ref{Table2} below.
The improvement with respect to the ranges found in \cite{HP} is particularly noticeable for the lower bounds
(recall that this is where the method of moments came into the picture).
\begin{table}[h]
\begin{center}
\begin{tabular}{|c||c||c|c||c|}
\hline
$\la(\cP_d)$ & known lower bound & our lower bound & our upper bound & known upper bound \\
\hline
$d=3$ & $1.3539$ & $1.35667...$ & $1.35696...$ & $1.3577$\\
\hline
$d=4$ & $1.4524$ & $1.45902$... & $1.45951...$ & $1.4611$\\
\hline	 
$d = 5$ & $1.525 $ & $1.53817...$ & $1.53895...$ & $1.543 $\\
\hline
$d=6$ & $1.580$ & $1.60271...$ & $1.60383...$ & $1.613$\\
\hline
$d=7$ & $1.624$ & $1.65693...$ &$1.65859...$  & $1.669$\\
\hline
$d=8$ & $1.660$ & $1.70483...$ & $1.70731...$ & $1.721$\\
\hline
$d=9$ & $1.678$ & $1.74774...$ & $1.75107...$ & $1.775$\\
\hline
$d=10$ & $1.696$ & $1.78658...$ & $1.79076...$ & $1.814$\\
\hline
$d=11$ & NA & $1.82169...$ & $1.82701...$ & NA\\
\hline
$d=12$ & NA & $1.85380...$ & $1.86216...$  & NA\\
\hline
\end{tabular}
\end{center}
\caption{For the spaces $\cP_d$ of polynomials of degree at most $d$,
lower and upper bounds on the projection constants obtained by our method compared to the ones stated in \cite{HP}.}
\label{Table2}
\end{table}
\end{samepage}

\begin{figure}[h]
\center
\subfigure{
\includegraphics[width=0.48\textwidth]{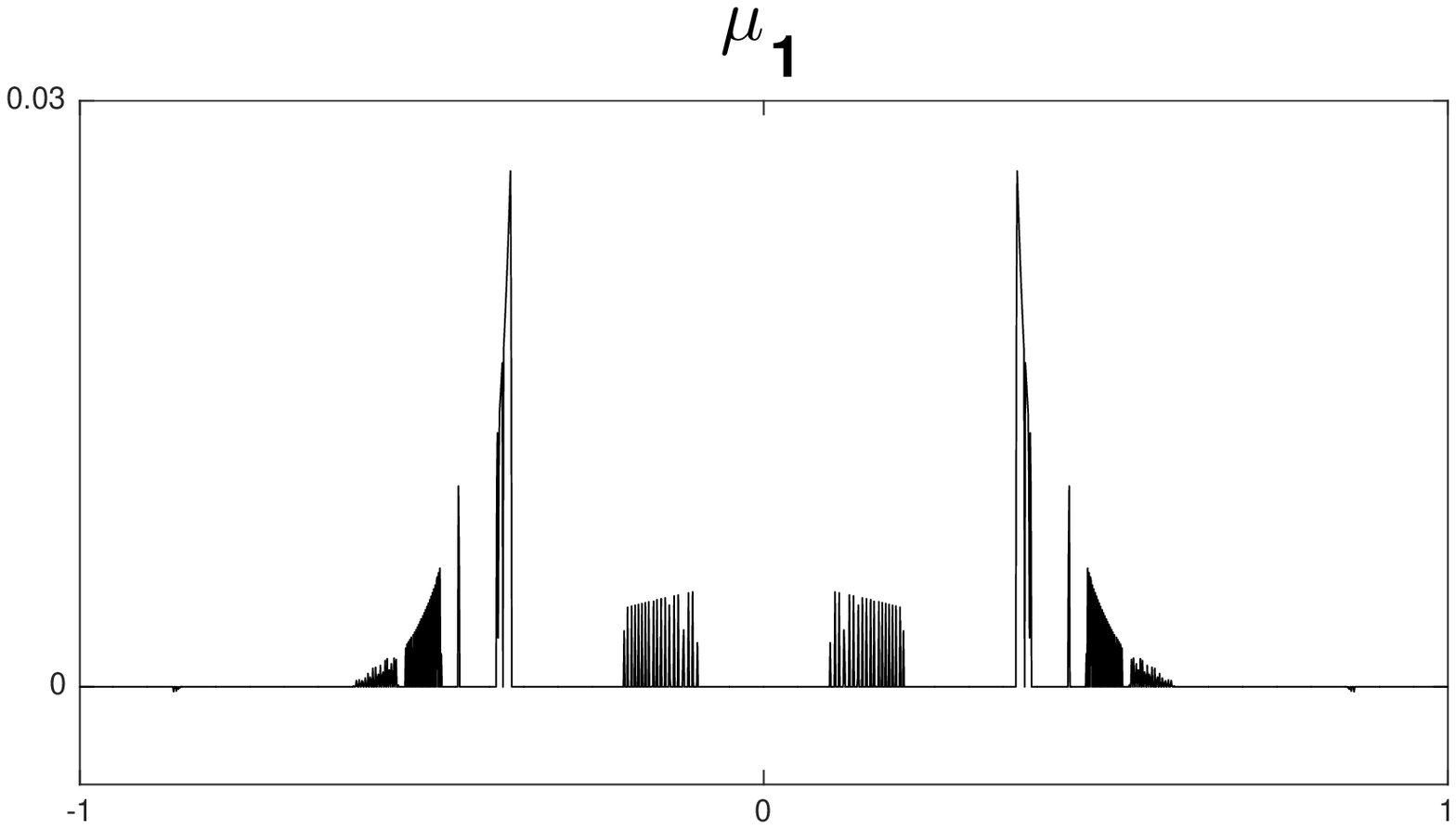}
}
\subfigure{
\includegraphics[width=0.48\textwidth]{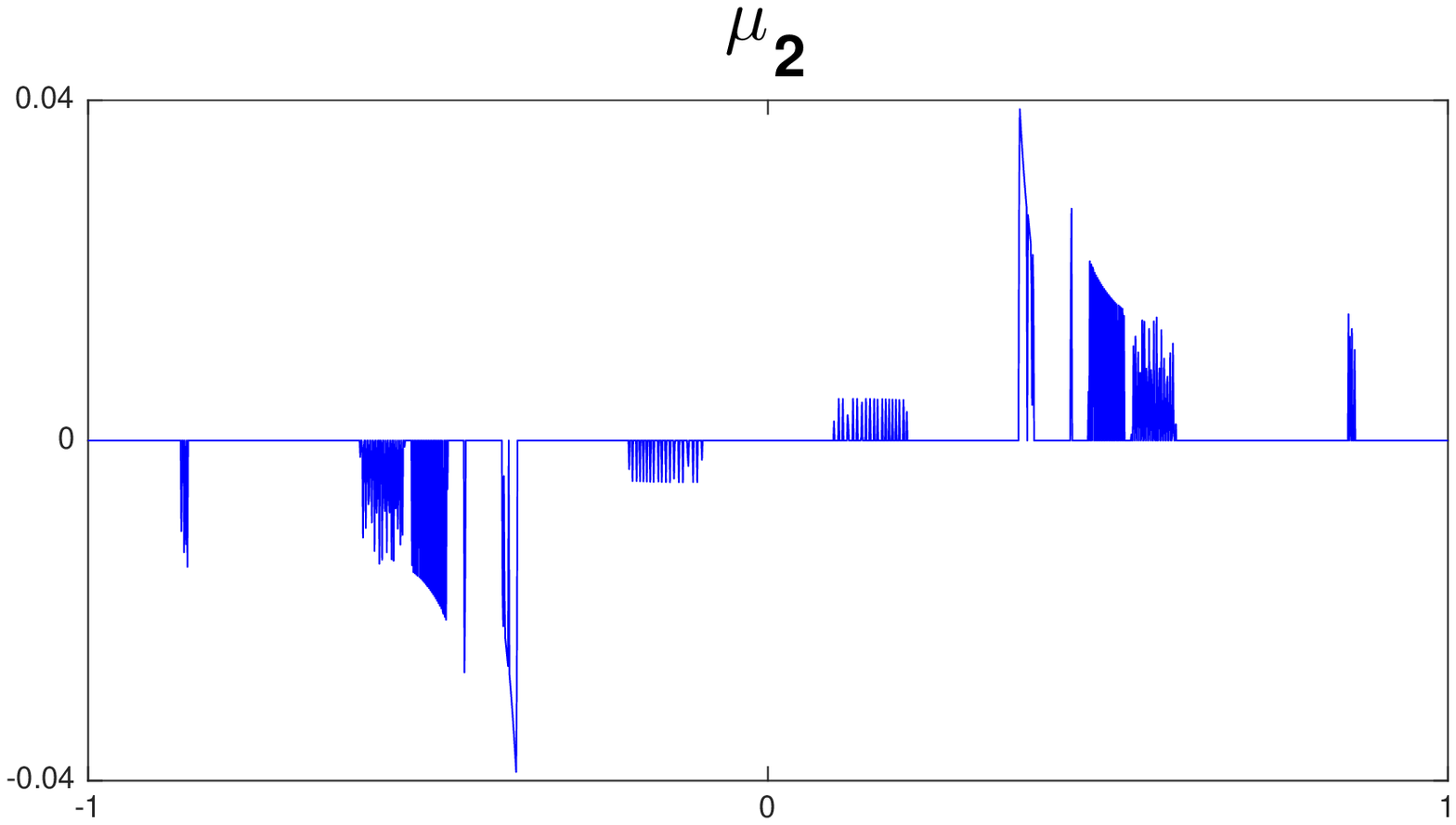}
}
\subfigure{
\includegraphics[width=0.48\textwidth]{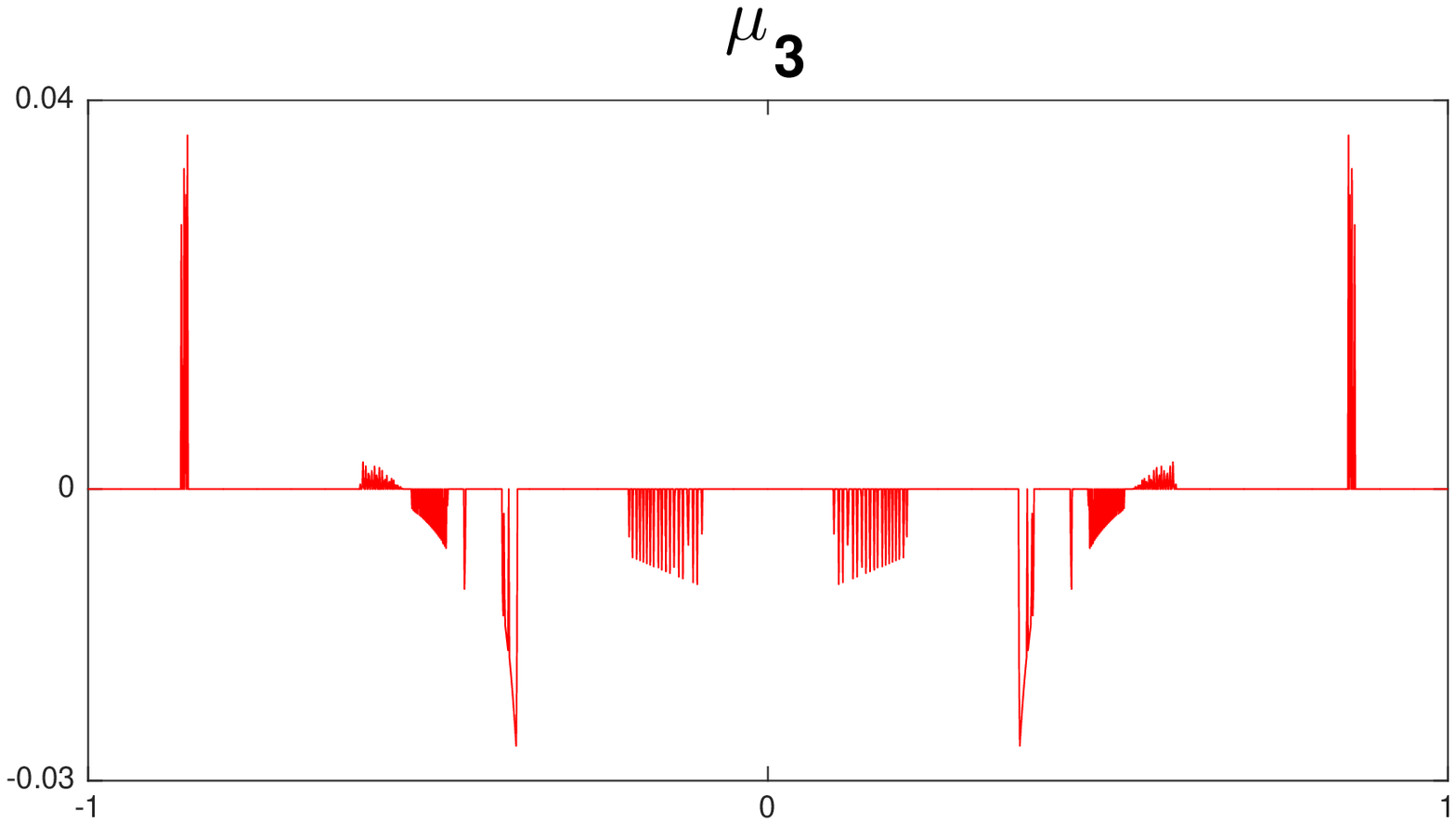}
}
\subfigure{
\includegraphics[width=0.48\textwidth]{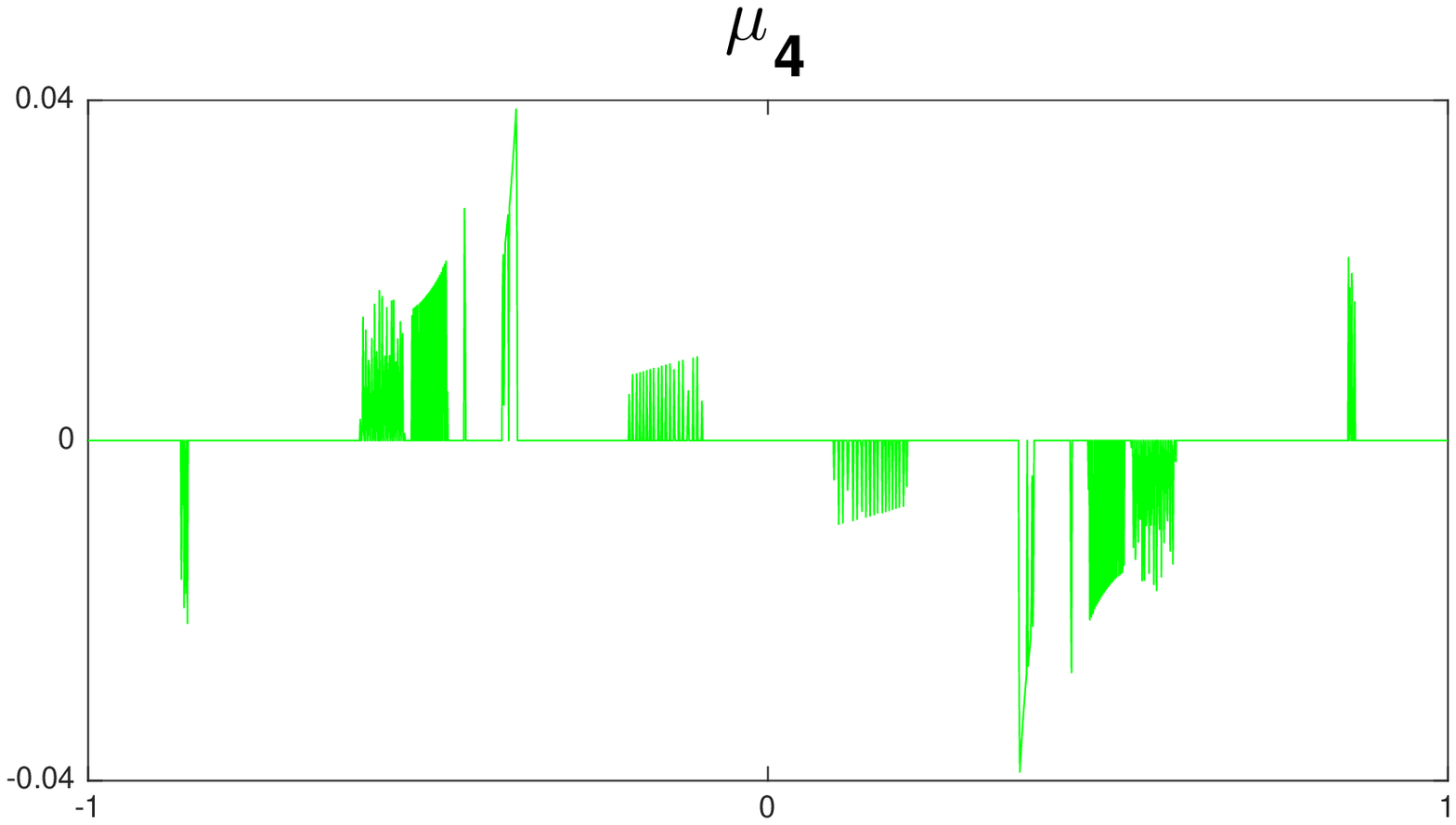}
}
\caption{The measures $\mu_1,\mu_2,\mu_3,\mu_4$
associated with the functionals dual to $1,x,x^2,x^3$
in the approximate minimal projection obtained by solving \eqref{ProgUBSym}.}
\label{Fig2}
\end{figure}

Let us now come to a close by examining the approximate minimal projection onto the cubics obtained by solving \eqref{ProgUBSym}.
For the measures $\mu_1,\mu_2,\mu_3,\mu_4$
associated with the functionals dual to $1,x,x^2,x^3$,
we detected atoms at $-1$ and $1$ (but none at $0$)
and Figure \ref{Fig2} indicates that their continuous parts seem to possess a common support strictly included in $[-1,1]$, as was the case for quadratics. 
But a disparity with the quadratics now occurs in terms of shape preservation.
It was conjectured in \cite{PCM}, and proved for $d=2$,
that (one of the) minimal projections onto $\cP_d$ preserve $d$-convexity.
This means that if $P$ is a minimal projection from $\cC[-1,1]$ onto $\cP_d$,
then,
for any $f \in \cC^d[-1,1]$,
\be
f^{(d)} \ge 0 \mbox{ on } [-1,1]
\overset{?}{\imp} (P(f))^{(d)} \ge 0 \mbox{ on } [-1,1],
\ee
or equivalently, writing $P(f)(x) = \sum_{m=1}^{d+1} \eta_m(f) x^{m-1}$,
\be
f^{(d)} \ge 0 \mbox{ on } [-1,1]
\overset{?}{\imp} \eta_{d+1}(f) \ge 0.
\ee
Our computations for $d=3$ give some insight that this conjecture should be false.
Indeed, 
if a minimal projection is approximated by our solution of \eqref{ProgUBSym},
which has the form
\be 
P(f)(x) = \sum_{m=1}^{4} \left( \sum_{k=-K}^K A_{m,k} f_k \right) x^{m-1},
\qquad f_k:= f\left( \f{k}{K} \right),
\ee 
and if the condition $f''' \ge 0$ on $[-1,1]$ is replaced by its discrete version 
\be
\Delta_3(f)_k := f_{k+3} - 3 f_{k+2} + 3 f_{k+1} - f_k
\ge 0,
\qquad k \in \ibt{-K}{K-3},
\ee
then, setting $a = A(4,:) \in \bR^{2K-1}$, the question becomes
\be
\min_{f \in \bR^{2K-1}} \left\{  \langle a, f \rangle 
 \; : \,
\Delta_3(f) \ge 0 \right\}
\; \overset{?}{\ge} 0.
\ee
This question is answered negatively
by solving a linear program.

\section{Outlook}

As a concluding message,
we reiterate our belief in the usefulness of modern optimization techniques for solving problems in Approximation Theory.
Purists will argue that `solving computationally' is not really solving,
but benefits are undeniable for building intuition about the problems at hand.
This article demonstrated, for instance,
how the method of moments elucidates the problem of minimal projections onto polynomial spaces
and it strongly hinted that minimal projections are not unique and do not preserve shape.
Our technique can be extended in several directions,
as long as linear-programming upper bounds and semidefinite-programming lower bounds match up to a desired accuracy.
For example, given subspaces $\cU \inc \cV$ of $\cC[-1,1]$,
one could find a linear map $P$ from  $\cC[-1,1]$ into $\cV$, not $\cU$,
with minimal norm among those satisfying $P(u) = u$ for all $u \in \cU$,
or even satisfying $P(u) = F(u)$ for all $u \in \cU$ with some $F \not= {\rm Id}_{\cU}$,
and shape-preservation properties may be added, etc...    
A particularly interesting situation concerns multivariate polynomial spaces,
for which nothing is known except the results of \cite{SS}.
Conceptually, the same technique applies,
but one quickly runs into numerical limitations.
One could throw in more computational power,
of course,
but it seems wiser to refine the method first,
possibly with a back-and-forth process between upper bound and lower bound.

\end{document}